\documentstyle{amsppt}
\magnification=1200
\catcode`\@=11
\redefine\logo@{}
\catcode`\@=13

\define \bn{\Bbb N}
\define \bz{\Bbb Z}
\define \bq{\Bbb Q}
\define \br{\Bbb R}
\define \bc{\Bbb C}

\define \M{{\Cal M}}

\define\La{{\Cal L}}
\define\geg{{\goth g}}
\define\0o{{\overline 0}}
\define\1o{{\overline 1}}
\define\rk{\text{rk}~}


\define\mult{\text{mult}}

\define\im{{\text{im}}}


\input epsf.tex

\TagsOnRight

\topmatter
\title
A Theory of Lorentzian Kac--Moody Algebras
\endtitle

\author
Viacheslav V. Nikulin \footnote{Supported by
Grant of Russian Fund of Fundamental Research
\hfill\hfill}
\endauthor

\address
Steklov Mathematical Institute,
ul. Gubkina 8, GSP-1 Moscow 117966 Russia
\endaddress

\email
slava\@nikulin.mian.su
\endemail

\dedicatory
Dedicated to 90-th Anniversary of Lev S. Pontryagin
\enddedicatory

\abstract
We present a variant of the Theory of Lorentzian (i. e. with
a hyperbolic generalized Cartan matrix) Kac--Moody algebras
recently developed by V. A. Gritsenko and the author. It is closely 
related with and strongly uses results of R. Borcherds. 
This theory should generalize 
well-known Theories of finite Kac--Moody
algebras (i. e. classical semi-simple Lie algebras corresponding to
positive generalized Cartan matrices) and affine Kac--Moody
algebras (corresponding to semi-positive generalized Cartan
matrices).

Main features of the Theory of Lorentzian Kac--Moody algebras are:
One should consider {\it generalized} Kac--Moody algebras 
introduced by Borcherds. Denominator
function should be an automorphic form on IV
type Hermitian symmetric domain (first example of this type related with
Leech lattice was found by Borcherds).
The Kac--Moody algebra is graded by elements of an integral
hyperbolic lattice $S$. Weyl group acts in the hyperbolic space related
with $S$ and has a fundamental
polyhedron $\Cal M$ of finite (or almost finite) volume and a  
lattice Weyl vector.

There are results and conjectures which permit (in principle) to get
a ``finite'' list of all possible Lorentzian Kac--Moody algebras.
Thus, this theory looks very similar to Theories of finite and
affine Kac--Moody algebras but is much more complicated.
There were obtained some classification 
results on Lorentzian Kac--Moody algebras and
many of them were constructed.
\endabstract

\rightheadtext
{Lorentzian Kac--Moody algebras}
\leftheadtext{V.V. Nikulin}

\endtopmatter

\document

\head
0. Introduction
\endhead

Lev Semenovich Pontryagin was the Great Mathematician. He had many
points of interest. One of them was topological groups,
Lie groups and Lie algebras. I mention his classical book 
``Topological Groups'' \cite{P}
containing his results.

On this conference, I am glad to present the talk related with
this subject. It is devoted to
{\it Lorentzian} (or {\it hyperbolic})
Kac--Moody Lie algebras which are a hyperbolic
analogy of classical semi-simple Lie algebras. A theory of
these algebras was recently developed by V.A. Gritsenko and the
author \cite{GN1}---\cite{GN7}.  
It is closely related with and strongly uses results
of R. Borcherds \cite{B1} --- \cite{B5}.

\head
1. Lorentzian Kac--Moody Algebras
\endhead

\subhead
1.1. Some general results on Kac--Moody algebras
\endsubhead
One can find all definitions and details in
the classical book by Victor Kac \cite{K}.

A {\it generalized Cartan matrix} $A$ is an integral square
matrix of a finite rank which has only $2$ on the diagonal and
non-positive integers out of the diagonal. 
We shall consider only {\it symmetrizable} 
generalized Cartan matrices $A$. It means that there exists a
diagonal matrix $D$ with positive rational coefficients 
such that $B=DA$ is integral and symmetric. Then $B$ is called the
{\it symmetrization} of $A$. By definition, 
$\text{sign}(A)=\text{sign}(B)$. We shall suppose that $A$ is
{\it indecomposable} which means that there does not exist
a decomposition $I=I_1\cup I_2$ of the set $I$ of indices of $A$
such that $a_{ij}=0$ if $i\in I_1$ and $j\in I_2$.

Each generalized Cartan matrix $A$ defines a Kac--Moody
Lie algebra $\geg(A)$ over $\bc$. 
The Kac--Moody algebra $\geg(A)$ is defined
by the set of generators and defining relations prescribed by the 
generalized Cartan matrix $A$. They are due to V. Kac and R. Moody. In
fact, they are a natural generalization of classical results of
Killing, Cartan, Weyl, Chevalley and Serre about finite-dimensional
semisimple Lie algebras. One should introduce
the set of generators $h_i$, $e_i$, $f_i$, $i \in I$,
with defining relations
$$
\cases
&[h_i,\,h_j]=0,\ \ [e_i,\,f_i]=h_i, \ \ [e_i,f_j]=0,\ \text{if\ }i\not=j,\\
&[h_i,\,e_j]=a_{ij}e_j,\ \ [h_i,\,f_j]=-a_{ij}f_j,\\
&(ad~e_i)^{1-a_{ij}}e_j=(ad~f_i)^{1-a_{ij}}f_j=0,\ \text{if}\ i\not=j.
\endcases
\tag{1.1}
$$
The important property of the algebra $\geg(A)$ is
that it is simple or almost simple: it is simple after factorization
by some known ideal.

We mention some general features of the theory of Kac--Moody
algebras $\geg(A)$.

1. The symmetrization $B$ defines a free $\bz$-module 
$Q=\sum_{i\in I}\bz\alpha_i$ with generators $\alpha_i$, $i\in I$, 
equipped with symmetric bilinear form $\left((\alpha_i,\,\alpha_j)\right)=B$ 
defined by the symmetrization $B$. 
The $Q$ is called {\it root lattice}. 
The algebra $\geg (A)$ is {\it graded by
the root lattice $Q$} (by definition, generators $h_i$, $e_i$, $f_i$ have
weights $0$, $\alpha_i$, $(-\alpha_i)$ respectively):
$$
\geg(A)=\bigoplus_{\alpha\in Q}{\geg_\alpha}=
\geg_0 \bigoplus \left(\bigoplus_
{\alpha\in \Delta_+}\geg_\alpha\right) \bigoplus
\left(\bigoplus_{\alpha\in \Delta_+}
\geg_{-\alpha}\right)
\tag{1.2}
$$
where $\geg_\alpha$ are finite dimensional linear spaces,
$[\geg_\alpha,\,\geg_\beta]\subset \geg_{\alpha +\beta}$,
$\geg_0\equiv Q\otimes \bc$ is commutative,
and is called {\it Cartan subalgebra}. 
An element $0\not=\alpha\in Q$ is called {\it root} if $\geg_\alpha\not=0$.
The $\geg_\alpha$ is called {\it root space} corresponding to $\alpha$. 
The dimension $\mult(\alpha)=\dim \geg_\alpha$ is called
{\it multiplicity} of the root $\alpha$. In \thetag{1.2}, 
$\Delta \subset Q$ is the set of all roots. It is divided in
the set of positive $\Delta_+\subset \sum_{i\in I}{\bz_+\alpha_i}$
and negative $-\Delta_+$ roots. A root $\alpha \in \Delta$ is called 
{\it real} if $(\alpha,\alpha)>0$. 
Otherwise, (if $(\alpha,\alpha)\le 0$), it is called {\it imaginary}.
Every real root $\alpha$ defines {\it a reflection}
$s_\alpha:x\mapsto x-\left(2(x,\alpha)/(\alpha,\alpha)\right)\alpha$, 
$x\in Q$. All refections $s_\alpha$ in real roots generate 
{\it Weyl group} $W\subset O(S)$. The set of roots $\Delta$ 
and multiplicities of roots are $W$-invariant. 

2. One has the {\it Weyl---Kac
denominator identity} which permits to
calculate multiplicities of roots: 
$$
e(-\rho)\prod_{\alpha\in \Delta_+}{(1-e(-\alpha))^{\mult (\alpha)}} =
\sum_{w\in W}{\det(w)e(-w(\rho))}.
\tag{1.3}
$$
Here $e(\,\cdot\, )\in \bz[Q]$ are formal exponents,   
$\rho$ is called {\it Weyl vector} and is defined by 
the condition $(\rho,\alpha_i)=-(\alpha_i,\alpha_i)/2$ for any $i\in I$. 

The formula \thetag{1.3} is combinatorial, and 
formulae for multiplicities $\mult(\alpha)$ are unknown in general. 
One approach to solve this problem is to replace the
formal function \thetag{1.3} by non-formal
one (e. g. replacing formal exponents by non-formal ones) 
to get a function with ``good'' properties. These good 
properties may help to find the formulae for multiplicities.

\subhead
1.2. Finite and affine cases
\endsubhead
There are two cases when we have very clear picture (or Theory)
of Kac--Moody algebras.

{\it Finite case:} The generalized Cartan matrix $A$ is positive definite,
$A>0$. Then $\geg(A)$ is finite-dimensional, and we get the classical
theory of finite-dimensional semisimple Lie algebras.

{\it Affine case:} The generalized Cartan matrix $A$ is semi-positive
definite, $A\ge 0$. Then $\geg(A)$ is called {\it affine}.

For both these cases we have three very nice properties:

(I). There exists the classification of all possible generalized
Cartan matrices $A$: They are classified by Dynkin (for finite case) 
and by extended Dynkin (for affine case) diagrams.

(II). In the denominator identity \thetag{1.2},  formal exponents may
be replaces by non-formal ones to give a function with nice properties:
For finite case this gives a polynomial. For affine case this gives a
Jacobi automorphic form. Using these properties (or directly), one
can find all multiplicities.

(III) Both these cases have extraordinary importance in
Mathematics and Phy\-sics.

\smallpagebreak

We want to construct similar Theory
(to the theories of finite and affine Lie
algebras above) for {\it Lorentzian (or hyperbolic) case} 
when the generalized Cartan matrix $A$ is {\it hyperbolic:
it has exactly one negative square, all its other squares are
either positive or zero.} There are plenty of hyperbolic
generalized Cartan matrices, it is impossible to
find all of them and classify. On the other hand, probably 
not all of them give interesting Kac--Moody algebras,
and one has to find natural conditions on these matrices.

\smallpagebreak

\subhead
1.3. Lorentzian case. Example of R. Borcherds
\endsubhead
We have the following {\it key example due to R. Borcherds} 
\cite{B1}---\cite{B5}.  

For Borcherds example, the root lattice $Q=S$ where $S$ is
an even hyperbolic unimodular lattice $S$ of signature
$(25,1)$. Here ``even'' means that $(x,x)$ is even for any $x\in S$.
``Unimodular'' means that the dual lattice $S^\ast$ coincides with  $S$,
equivalently, for a bases $e_1,\dots, e_{26}$ of $S$
the determinant of the Gram matrix $\big((e_i,\,e_j)\big)$
is equal to $\pm 1$. A lattice $S$ with these properties
is unique up to isomorphism.
For Borcherds example, the Weyl group $W$ is generated by reflections
$s_\alpha:x\mapsto x-(x,\alpha)\alpha$, $x\in S$, in all elements 
$\alpha \in S$ with $\alpha^2=2$.  
The group $W$ is discrete in the hyperbolic space
$\La (S)=V^+(S)/\br_{++}$. Here $V^+(S)$ is a half of the light cone
$V(S)=\{x\in S\otimes \br \ |\ x^2<0\}$ of $S$.  
The $\La (S)$ is the set of rays in $V^+(S)$. 

A fundamental chamber $\M\subset \La (S)$ for $W$ is defined by the set $P$
of elements $\alpha\in S$ with $\alpha^2=2$ which are {\it orthogonal to}  
$\M$. It has the following description 
due to J. Conway \cite{Co}. 
There exists an orthogonal
decomposition $S=[\rho,\,e]\oplus L$ where the Gram matrix
of elements $\rho,\,e$ is equal to
$U=\left(\matrix0&-1\\-1&0\endmatrix\right)$
(in particular, $(\rho,\rho)=0$),
and $L$ is the Leech lattice, i. e. positive definite even unimodular
lattice of the rank $24$ without elements with square $2$. 
The set $P$ of roots
which are orthogonal to the fundamental chamber $\M$ (or the
set of {\it simple roots}) of $W$ is equal to 
$$
P=\{\alpha \in S\ |\ (\alpha,\,\alpha)=2\ \and\ (\rho,\alpha)=-1 \}.
\tag{1.4}
$$
It means that the fundamental chamber $\M\subset \La(S)$ is equal to
$$
\M=\{\br_{++}x\in\La(S)\ |\ (x,\ P)\le 0\} 
\tag{1.5}
$$
and $P$ is a minimal set with this property. We mention that the
fundamental polyhedron $\M$ has {\it "almost finite" volume}. It means
that $\M$ is finite in any angle with the center at infinity
$\br_{++} \rho$ of the hyperbolic space $\La(S)$.

The matrix
$$
A=\bigl((\alpha,\,\alpha^\prime )\bigr),\ \ \alpha,\,\alpha^\prime \in P  
\tag{1.6}
$$
is a generalized Cartan matrix and $\rho$ is the Weyl vector:
$$
(\rho,\alpha)=-(\alpha,\alpha)/2, \ \forall \alpha \in P.
\tag{1.7}
$$
We have the classical $SL_2(\bz)$-modular cusp
form $\Delta$ of the weight $12$ on the upper-half plane $\text{im}~q>0$:
$$
\Delta=q\prod_{n=1}^{\infty}(1-q^n)^{24}=
\sum_{m\ge 0}\tau(m)q^m,
\tag{1.8}
$$
where $q=exp(2\pi i \tau)$. We have
$$
\Delta^{-1}=\sum_{n\ge 0}p_{24}(n)q^{n-1}
\tag{1.9}
$$
where $p_{24}(n)$ are positive integers.
Borcherds \cite{B2} proved the {\it identity}
$$
\Phi(z)=\exp{(-2\pi i (\rho,z))}\prod_{\alpha\in \Delta_+}
{(1-\exp{(-2\pi i (\alpha,z))})^{p_{24}(1-(\alpha,\alpha)/2)}}=
$$
$$
\sum_{w\in W}{\det(w)\sum_{m>0}{\tau (m)\exp{(-2\pi i (w(m\rho),z))}}}.
\tag{1.10}
$$
Here $\Delta_+=\{\alpha \in S\,|\,\alpha^2=2\ \and\ (\alpha,\,\rho)<0\}\cup
(S\cap \overline{V^+(S)}-\{0\})$. The variable
$z$ runs through the {\it complexified cone}
$\Omega(V^+(S))=S\otimes\br\,+\,i V^+(S )$ of the light cone
$V^+(S)$.  Moreover, Borcherds \cite{B4}, \cite{B5} proved that the
function $\Phi(z)$ {\it is an automorphic form of weight $12$}
with respect to the
group $O^+(T)$ where $T=U\oplus S$ is the extended lattice of
the signature $(26,2)$. The group $O^+(T)$ naturally acts in the
Hermitian symmetric domain of type IV 
$$
\Omega(T)=\{\bc\omega \subset T\otimes \bc\ |\ (\omega,\,\omega)=0\ \and\
(\omega,\,\overline{\omega})<0 \}_0,
\tag{1.11}
$$
which has canonical identification with $\Omega(V^+(S))$ as follows:
$z\in\Omega (V^+(S))$ defines the element 
$\bc\omega_z\in \Omega(T)_0$
where $\omega_z =\big((z,z)/2\big)e_1+
e_2\oplus z \in T\otimes \bc$
and $e_1$, $e_2$ is the bases of the lattice $U$ with
the Gram matrix $U$ above. Here ``automorphic of
the weight 12'' means that the function
$\widetilde{\Phi}(\lambda \omega_z)=\lambda^{-12}\Phi(z)$, 
$\lambda \in \bc^\ast$, 
is homogeneous of the degree $-12$ (it is obvious)
in the cone
$\widetilde{\Omega(T)}_0$ over $\Omega(T)_0$,
and $\widetilde{\Phi}(g\omega)=\widetilde{\Phi}(\omega)$
for any $g\in O^+(T)$ where $O^+(T)$ is the subgroup of index $2$ of
the group $O(T)$ which keeps the connected component \thetag{1.11}
(marked by $0$).

The identity \thetag{1.10} looks very familiar to
the form \thetag{1.3} of the denominator identity for Kac--Moody
algebras, but it has some difference.

To interpret \thetag{1.10} as a denominator identity of a Lie
algebra, Borcherds introduced \cite{B1} 
{\it generalized Kac--Moody algebras}  
$\geg(A^\prime )$ which correspond to more general matrices $A^\prime$ than
generalized Cartan matrices. Here I will call them as 
{\it generalized generalized Cartan matrices}.
Difference is that a generalized
generalized Cartan matrix $A^\prime$ may also have 
non-positive real elements 
$a_{ij}\le 0$ on the diagonal and out of the diagonal, 
but all $a_{ij}\in \bz$ if $a_{ii}=2$. 
A definition of the generalized Kac--Moody 
algebra $\geg (A^\prime )$ corresponding to a generalized 
generalized Cartan matrix $A^\prime$ is similar to \thetag{1.1}. 
One should replace the last line of \thetag{1.1} by
$$
(ad~e_i)^{1-a_{ij}}e_j=(ad~f_i)^{1-a_{ij}}f_j=0\  
\text{if}\ i\not=j\ \text{and}\ a_{ii}=2,
\tag{1.12}
$$
and add the relation
$$
[e_i,\,e_j]=[f_i,\,f_j]=0\ \text{if}\ a_{ij}=0.
\tag{1.13}
$$
Borcherds showed that generalized Kac--Moody algebras have similar 
properties to ordinary Kac--Moody algebras.
They also have a
denominator identity which has more general form than \thetag{1.3}
and includes \thetag{1.10} as a particular case.

The identity \thetag{1.10} is the denominator identity for
the generalized Kac--Moody algebra $\geg(A')$ where $A'$ is the generalized
generalized Cartan matrix equals to the Gram
matrix $A^\prime=\bigl((\alpha,\,\alpha^\prime )\bigr)$, 
$\alpha,\,\alpha^\prime \in P^\prime$ where  
$$
P^\prime=P\cup 24\rho \cup 24(2\rho) \cup \dots \cup 24(n\rho)\cup\dots
\tag{1.14}
$$
is the sequence of elements of the lattice $S$.  
Here $24(n\rho)$ means that we take the element $n\rho$ twenty four
times to get the Gram matrix $A^\prime$. See details in
\cite{B1}---\cite{B3}.   

In \thetag{1.14}, the set $P^\prime$ defining $A^\prime$
is called the {\it set of simple roots}. It is divided in the set
${P^\prime}^{re}=P$ of {\it simple real roots} (they are orthogonal
to the fundamental chamber $\M$ of the Weyl group $W$ and have 
positive square) and is the same
as for the ordinary Kac--Moody algebra $\geg(A)$ defined by the
generalized Cartan matrix $A$ in \thetag{1.6}. The additional sequence 
$$
{P^\prime}^{im}=24\rho \cup 24(2\rho) \cup \dots \cup 24(n\rho)\cup\dots
\tag{1.15}
$$
of $P^\prime$ (elements of ${P^\prime}^{im}$ have non-positive square) 
is defined by the Fourier coefficients in the sum part of
the identity \thetag{1.10}. For example, $24$ is defined by the $24$ in
\thetag{1.8}. Together ${P^\prime}^{re}$ and ${P^\prime}^{\im}$
define the generalized generalized Cartan matrix
$A^\prime$ and the generalized Kac--Moody algebra $\geg(A^\prime)$.

Borcherds example is very fundamental and beautiful. It has several
applications in Mathematics (e.g. for the Monster) and Physics
(e. g. for String Theory).

\subhead
1.4. A Theory of Lorentzian Kac--Moody algebras
\endsubhead
Analizing Borcherds example, one can suggest a general
class of Lorentzian Kac--Moody algebras (or automorphic Kac--Moody
algebras) $\geg$, see \cite{N7}, \cite{N8}, \cite{GN1}---\cite{GN7}. 

One takes {\it data (1)---(5)} below: 

(1) A {\it hyperbolic lattice} $S$ (i. e. an integral symmetric bilinear
form of signature $(n,1)$).

(2) A {\it reflection group} (or {\it Weyl group})   
$W\subset O(S)$ generated by reflections in 
roots of $S$. We remind that $\alpha\in S$ is called {\it root} if
$\alpha^2>0$ and $\alpha^2\,|\,2(\alpha,S)$. Any root defines
a reflection
$s_\alpha: x\mapsto x-(2(x,\alpha)/\alpha^2)\alpha$, $x\in S$, which
gives an automorphism of the lattice $S$.

(3) {\it A set $P$ of orthogonal roots} to a {\it fundamental chamber}
$\M\subset \La(S)$ of $W$. It means that the set $P$ 
of roots of $S$ should have the property \thetag{1.5} 
and should be minimal
having this property. Moreover, the set $P$ should have a  
{\it Weyl vector} $\rho\in S\otimes \bq$
defined by the equality \thetag{1.7}
(it is more right to call it as {\it lattice Weyl vector}).

The main invariant of the data (1)---(3) is the {\it generalized 
Cartan matrix} 
$$
A=\left({2(\alpha, \alpha^\prime)\over (\alpha,\alpha)}\right),\ \ \
\alpha,\,\alpha^\prime\in P.
\tag{1.16}
$$
It defines data (1)---(3) up to some very clear equivalence and
defines the set of real roots of the algebra $\geg$ we want to construct.

(4) An {\it automorphic (holomorphic)
form} $\Phi(z)$ on a IV type Hermitian symmetric
domain, $z\in \Omega(V^+(S))=\Omega (T)$, with respect to a
subgroup $G\subset O^+(T)$ of finite index of an extended
lattice $T=U(k)\oplus S$ where
$U(k)=\left(\matrix0&-k\\-k&0\endmatrix\right)$, $k\in \bn$. 
(See \cite{GN5} for more general definition.) It
should have Fourier expansion of the form
of the denominator identity for a generalized Kac--Moody algebra
with hyperbolic generalized generalized Cartan matrix. This form is
$$
\Phi(z)=\sum_{w\in W}{\det(w)
\Bigl(\exp{\left(-2\pi i (w(\rho),z)\right)}\ -\hskip-10pt
\bigr.\sum_{a\in S\cap \br_{++}\M}
{\bigl.m(a)\exp{\left(-2\pi i (w(\rho+a),z)\right)}\Bigr)}}
\tag{1.17}
$$
where all coefficients $m(a)$ should be integral. The automorphic
form $\Phi$ defines the set of imaginary roots of the algebra $\geg$.

Like for Borcherds example, {\it data (1)---(4) define a generalized
Kac--Moody algebra or superalgebra (if some Fourier
coefficients $m(a)$ are negative)} $\geg$.  
See the definition of $\geg$ below. 
Using automorphic properties
of $\Phi(z)$, it is good to {\it calculate the product part of
the denominator identity}
$$
\Phi(z)=\exp{(-2\pi i (\rho,z))}
\prod_{\alpha\in \Delta_+}{\Bigl(1-\exp{\left(-2\pi i (\alpha,z)\right)}
\Bigr)^{\mult(\alpha)}} 
\tag{1.18}
$$
which gives multiplicities $\mult(\alpha)$ of roots $\alpha$ of $\geg$. 
For superalgebras case, the 
$\mult(\alpha)=\dim \geg_{\alpha,\0o}-\dim \geg_{\alpha,\1o}$ is
the difference of dimensions of even and odd parts of the root
space $\geg_\alpha$. 

It was understood that it is  
natural to suppose (at least, to have finiteness results)
the additional condition: 

(5) The automorphic form $\Phi$ on the domain $\Omega (V^+(S))=\Omega (T)$
should be {\it reflective}. It means that the divisor (of zeros) of $\Phi$
is union of quadratic divisors orthogonal to roots of
the extended lattice $T$. Here for a root
$\alpha \in T$ (the definition of a root of $T$ is the same as for
the lattice $S$) the {\it quadratic divisor orthogonal to} $\alpha$ is equal to
$$
D_\alpha=\{\bc\omega \in \Omega (T)\ |\ (\omega,\,\alpha)=0\}.
\tag{1.19}
$$

The property (5) is valid for Borcherds example above and in all known cases.
Moreover, it is true in the neighbourhood of the cusp where
the product \thetag{1.18} convergers. Thus, we want it to be true
globally.

A generalized Kac--Moody superalgebra $\geg$ corresponding to 
data (1)---(4) is given by the sequence $P^\prime \subset S$ 
{\it of simple roots}. 
This sequence is 
divided in a set ${P^\prime}^{re}$ 
of {\it simple real roots} and a sequence ${P^\prime}^{im}$ 
of {\it simple imaginary roots}. The sequence 
${P^\prime}^{im}$ is divided in the sequence  
${P^\prime}^{im}_{\0o}$ of {\it even simple imaginary roots} and a sequence  
${P^\prime}^{im}_{\1o}$ 
of {\it odd simple imaginary roots}. 
For any primitive $a\in S\cap \br_{++}\M$ 
with $(a,\,a)=0$, one should find $\tau(na)\in \bz$, $n\in \bn$, 
from the identity with the formal variable $T$:
$$
1-\sum_{k\in \bn}{m(ka)T^k}=\prod_{n\in \bn}{(1-T^n)^{\tau(na)}}.
$$
We set ${P^\prime}^{re}=P$ where $P$ is defined in (3). The set 
${P^\prime}^{re}$ is considered to be even: 
${P^\prime}^{re}={P^\prime}^{re}_{\0o}$ and 
${P^\prime}^{re}_{\1o}=\emptyset$. We set  
$$
\split 
{P^\prime}^{im}_{\0o}&=\{m(a)a\ |\  
a\in S\cap\br_{++}\M,\ (a,a)<0\ \text{and}\   
m(a)>0\}\cup\\
&\cup \{\tau(a)a\ |\ a\in S\cap\br_{++}\M,\ (a,a)=0\ \text{and}\  
\tau (a)>0\};
\endsplit
$$
$$
\split
{P^\prime}^{im}_{\1o}&=
\{-m(a)a\ |\ a\in S\cap\br_{++}\M,\ (a,a)<0\ \text{and}\ 
m(a)<0\}\cup\\
&\cup \{\tau(a)a\ |\ a\in S\cap\br_{++}\M,\ (a,a)=0\ \text{and}\ 
\tau (a)<0\}.
\endsplit 
$$
The generalized Kac--Moody superalgebra $\geg$ is a Lie superalgebra 
generated by $h_r$, $e_r$, $f_r$ where $r \in P^\prime$. 
All generators  $h_r$ are even, generators $e_r$, $f_r$ are even 
(respectively odd) if $r$ is even (respectively odd). They have the  
defining relations 1) --- 5) below:

\smallpagebreak  
 
1) The map $r\mapsto h_r$ for $r\in P^\prime$ gives an embedding of 
$S\otimes \bc$ to $\geg$ as an abelian subalgebra (it is even).  

2) $[h_r,\,e_{r^\prime}]=(r,\,r^\prime)e_{r^\prime}$ and 
$[h_r,f_{r^\prime}]=-(r,\,r^\prime)f_{r^\prime}$. 

3) $[e_r,\,f_{r^\prime}]=h_r$ if $r=r^\prime$, and is $0$ if 
$r\not=r^\prime$.
 
4) $(ad~e_r)^{1-2(r,\,r^\prime)/(r,\,r)}e_{r^\prime}=
(ad~f_r)^{1-2(r,\,r^\prime)/(r,\,r)}f_{r^\prime}=0\ 
\text{if $r\not=r^\prime$ and $(r,r)>0$}$\hfil\hfil
\newline
\phantom{(4)(4)} (equivalently, $r\in {P^\prime}^{re}$).

5) If $(r,\,r^\prime)=0$, then $[e_r,\,e_{r^\prime}]=
[f_r,\,f_{r^\prime}]=0$. 

See \cite{B1}, \cite{GN1}, \cite{GN2}, \cite{GN5}, \cite{R} 
for details.  We remark that this definition is equivalent 
(for ordinary generalized Kac--Moody algebras) to the definition above  
using a generalized generalized Cartan matrix defined by the Gram 
matrix of the sequence $P^\prime$.  
\smallpagebreak

{\it
The generalized Kac--Moody superalgebras $\geg$ above 
given by the data (1)---(5) constitute the Theory
of Lorentzian Kac--Moody algebras (or automorphic Lorentzian Kac--Moody
algebras) which we consider.} 

\smallpagebreak

By (4), they have similar property to the property (II) for 
finite and affine algebras: their denominator identity gives
an automorphic form. For Lorentzian case, 
it is an automorphic form on IV type symmetric domain.

What about a similar property to the property (I) 
for finite and affine algebras? How many data (1)---(5) one may have? 

Further, we suppose that $\rk S\ge 3$.
This condition is equivalent to considering
hyperbolic analogues of non-trivial finite-dimensional semi-simple
Lie algebras. When $\rk S =1,\ 2$, classification problem of data
(1)---(5) is different and, it seems, more simple. 

We have (see \cite{N3}, \cite{N4}, \cite{N7}, \cite{N8}, \cite{GN5}):   

\proclaim{Theorem 1} If $\rk S\ge 3$, the set of possible data (1)---(3)
in data (1)---(4) is finite when $(\rho,\rho)<0$
and is in essential finite when
$(\rho,\rho)=0$. The inequality $(\rho,\rho)>0$ is impossible.
\endproclaim

In the theorem, ``in essential finite'' 
means that the set might be infinite, 
but we have very clear understanding of the set. For example, the set of
possible Dynkin diagrams of the type $A_n$ is infinite, but we understand
the set very clearly.

Key point of the proof of Theorem 1 is that (1)---(4) imply
that $(\rho,\,\rho)\le 0$ and
the fundamental chamber $\M$ has finite (if $(\rho,\,\rho)<0$)
or almost finite (if $(\rho,\,\rho)=0$) volume, see \cite{N8}, \cite{GN5}. 
(Here ``almost finite'' means the same as for the Borcherds example above.)
Then the number of possible root lattices $S$ is finite. It
follows from old results of author \cite{N3}, \cite{N4} (and also
\cite{N8}) and \'E.B. Vinberg \cite{V}.
If there additionally exists a Weyl vector $\rho$ for $P$,
one also has finiteness (when $(\rho,\rho)<0$) or in
essential finiteness (when $(\rho,\rho)=0$) of the
sets of Weyl groups $W$, fundamental chambers
$\M$ (up to action of $W$) and the sets of orthogonal roots
$P$ to $\M$, see \cite{N8}. 
This gives finiteness or in essential finiteness
of possible generalized Cartan
matrices $A$ in \thetag{1.16}
corresponding to systems of simple real roots. 

It follows that in principle we can classify
all possible data (1)---(3) in the data (1)---(4).
It makes the Theory of Lorentzian Kac--Moody algebras 
similar to the Theories of finite and affine algebras.

It would be nice to have also finiteness results for the data (4), (5).
Recently, here,
there were obtained some partial finiteness results \cite{N9}, \cite{GN5}
which show that automorphic forms $\Phi(z)$ in (4) and (5) are
extremely rare. It makes very possible the following statement: 

\proclaim{Conjecture 2} If $\rk S\ge 3$,
the set of possible data (4), (5) is in essential finite.
\endproclaim

The reason why we expect the statement of Conjecture 2, 
is based on the {\it Koecher principle} (e. g. see \cite{Ba}):
Any holomorphic automorphic form on a Hermitian symmetric
domain $\Omega$ should have zeros in $\Omega$ if
$\dim \Omega - \dim \Omega_\infty\ge 2$.

Applying this principle to restrictions $\Phi|\Omega(T_1)$ of a
reflective automorphic form $\Phi$ on all subdomains
$\Omega(T_1)\subset \Omega (T)$ where $T_1\subset T$ is a sublattice
of $T$ of signature $(k,2)$,
one obtains very strong restrictions on the lattice $T$ to have a
reflective automorphic form $\Phi$.
This was shown in \cite{N9}. 

We think that Conjecture 2 is very interesting. From our point of
view, the theory of reflective automorphic forms on IV type domains
$\Omega (T)$ where $T$ is a lattice of signature $(n,2)$ is
similar (Mirror Symmetric)
to the theory of reflection groups $W$ of hyperbolic lattices
$S$ with fundamental
chamber of finite or almost finite volume, see \cite{GN3}, 
\cite{GN5}---\cite{GN7}, \cite{N9}. 

It is interesting to classify (or describe)
the conjecturally ``finite set'' 
of data (1)---(5). Even a finite set may have very interesting structure.
As a result, we will have a Theory of Lorentzian
Kac--Moody algebras which one can consider as a hyperbolic analogy of
the Theories of finite and affine Kac--Moody algebras.

\smallpagebreak

At the end, we describe a small piece of this classification
which was obtained in \cite{GN4}---\cite{GN6}.

There are exactly 12 generalized Cartan matrices $A$ of
data (1)---(3) in (1)---(4)
which are symmetric, have the rank 3 and have the
Weyl vector $\rho$ with $(\rho,\rho)<0$. They are given below:

\smallpagebreak

\centerline{\bf The list of all symmetric hyperbolic generalized
Cartan matrices}
\centerline{\bf of the rank $3$ with   
$vol(\M)<\infty$ and a lattice Weyl vector
$\rho$}

$$
A_{1,0}=
\left(\smallmatrix
\hphantom{-}{2}&\hphantom{-}{0}&{-1}\cr
\hphantom{-}{0}&\hphantom{-}{2}&{-2}\cr
{-1}&{-2}&\hphantom{-}{2}\cr
\endsmallmatrix\right),\ \
A_{1,I}=
\left(\smallmatrix
\hphantom{-}{2}&{-2}&{-1}\cr
{-2}&\hphantom{-}{2}&{-1}\cr
{-1}&{-1}&\hphantom{-}{2}\cr
\endsmallmatrix\right),\ \
A_{1,II}=
\left(\smallmatrix
\hphantom{-}{2}&{-2}&{-2}\cr
{-2}&\hphantom{-}{2}&{-2}\cr
{-2}&{-2}&\hphantom{-}{2}\cr
\endsmallmatrix\right),
$$
$$
A_{1,III}=
\left(\smallmatrix
\hphantom{-}{2}&{-2}&{-6}&{-6}&{-2}\cr
{-2}&\hphantom{-}{2}&\hphantom{-}{0}&{-6}&{-7}\cr
{-6}&\hphantom{-}{0}&\hphantom{-}{2}&{-2}&{-6}\cr
{-6}&{-6}&{-2}&\hphantom{-}{2}&{0}\cr
{-2}&{-7}&{-6}&\hphantom{-}{0}&\hphantom{-}{2}\cr
\endsmallmatrix\right);\ \
A_{2,0}=
\left(\smallmatrix
\hphantom{-}{2}&{-2}&{-2}\cr
{-2}&\hphantom{-}{2}&\hphantom{-}{0}\cr
{-2}&\hphantom{-}{0}&\hphantom{-}{2}\cr
\endsmallmatrix\right),\ \
A_{2,I}=
\left(\smallmatrix
\hphantom{-}{2}&{-2}&{-4}&\hphantom{-}{0}\cr
{-2}&\hphantom{-}{2}&\hphantom{-}{0}&{-4}\cr
{-4}&\hphantom{-}{0}&\hphantom{-}{2}&{-2}\cr
\hphantom{-}{0}&{-4}&{-2}&\hphantom{-}{2}\cr
\endsmallmatrix\right),
$$
$$
A_{2,II}=
\left(\smallmatrix
\hphantom{-}{2}&{-2}&{-6}&{-2}\cr
{-2}&\hphantom{-}{2}&{-2}&{-6}\cr
{-6}&{-2}&\hphantom{-}{2}&{-2}\cr
{-2}&{-6}&{-2}&\hphantom{-}{2}\cr
\endsmallmatrix\right),\ \
A_{2,III}=
\left(\smallmatrix
\hphantom{-}{2}&{-2}&{-8}&{-16}&{-18}&{-14}&{-8}&\hphantom{-}{0}\cr
{-2}&\hphantom{-}{2}&\hphantom{-}{0}&{-8}&{-14}&{-18}&{-16}&{-8}\cr
{-8}&\hphantom{-}{0}&\hphantom{-}{2}&{-2}&{-8}&{-16}&{-18}&{-14}\cr
{-16}&{-8}&{-2}&\hphantom{-}{2}&\hphantom{-}{0}&{-8}&{-14}&{-18}\cr
{-18}&{-14}&{-8}&\hphantom{-}{0}&\hphantom{-}{2}&{-2}&{-8}&{-16}\cr
{-14}&{-18}&{-16}&{-8}&{-2}&\hphantom{-}{2}&\hphantom{-}{0}&{-8}\cr
{-8}&{-16}&{-18}&{-14}&{-8}&\hphantom{-}{0}&\hphantom{-}{2}&{-2}\cr
\hphantom{-}{0}&{-8}&{-14}&{-18}&{-16}&{-8}&{-2}&\hphantom{-}{2}\cr
\endsmallmatrix\right);
$$
$$
A_{3,0}=
\left(\smallmatrix
\hphantom{-}{2}&{-2}&{-2}\cr
{-2}&\hphantom{-}{2}&{-1}\cr
{-2}&{-1}&\hphantom{-}{2}\cr
\endsmallmatrix\right),\ \
A_{3,I}=
\left(\smallmatrix
\hphantom{-}{2}&{-2}&{-5}&{-1}\cr
{-2}&\hphantom{-}{2}&{-1}&{-5}\cr
{-5}&{-1}&\hphantom{-}{2}&{-2}\cr
{-1}&{-5}&{-2}&\hphantom{-}{2}\cr
\endsmallmatrix\right),
$$
$$
A_{3,II}=
\left(\smallmatrix
\hphantom{-}{2}&{-2}&{-10}&{-14}&{-10}&{-2}\cr
{-2}&\hphantom{-}{2}&{-2}&{-10}&{-14}&{-10}\cr
{-10}&{-2}&\hphantom{-}{2}&{-2}&{-10}&{-14}\cr
{-14}&{-10}&{-2}&\hphantom{-}{2}&{-2}&{-10}\cr
{-10}&{-14}&{-10}&{-2}&\hphantom{-}{2}&{-2}\cr
{-2}&{-10}&{-14}&{-10}&{-2}&\hphantom{-}{2}\cr
\endsmallmatrix\right),
$$
$$
A_{3,III}=
\left(\smallmatrix
\hphantom{-}{2}&{-2}&{-11}&{-25}&{-37}&{-47}&{-50}
&{-46}&{-37}&{-23}&{-11}&{-1}\cr
{-2}&\hphantom{-}{2}&{-1}&{-11}&{-23}&{-37}&
{-46}&{-50}&{-47}&{-37}&{-25}&{-11}\cr
{-11}&{-1}&\hphantom{-}{2}&{-2}&{-11}&{-25}
&{-37}&{-47}&{-50}&{-46}&{-37}&{-23}\cr
{-25}&{-11}&{-2}&\hphantom{-}{2}&{-1}&{-11}
&{-23}&{-37}&{-46}&{-50}&{-47}&{-37}\cr
{-37}&{-23}&{-11}&{-1}&\hphantom{-}{2}&{-2}
&{-11}&{-25}&{-37}&{-47}&{-50}&{-46}\cr
{-47}&{-37}&{-25}&{-11}&{-2}&\hphantom{-}{2}
&{-1}&{-11}&{-23}&{-37}&{-46}&{-50}\cr
{-50}&{-46}&{-37}&{-23}&{-11}&{-1}&\hphantom{-}{2}
&{-2}&{-11}&{-25}&{-37}&{-47}\cr
{-46}&{-50}&{-47}&{-37}&{-25}&{-11}&{-2}
&\hphantom{-}{2}&{-1}&{-11}&{-23}&{-37}\cr
{-37}&{-47}&{-50}&{-46}&{-37}&{-23}&{-11}&{-1}
&\hphantom{-}{2}&{-2}&{-11}&{-25}\cr
{-23}&{-37}&{-46}&{-50}&{-47}&{-37}
&{-25}&{-11}&{-2}&\hphantom{-}{2}&{-1}&{-11}\cr
{-11}&{-25}&{-37}&{-47}&{-50}&{-46}&{-37}
&{-23}&{-11}&{-1}&\hphantom{-}{2}&{-2}\cr
{-1}&{-11}&{-23}&{-37}&{-46}&{-50}
&{-47}&{-37}&{-25}&{-11}&{-2}&\hphantom{-}{2}\cr
\endsmallmatrix\right).
$$

\smallpagebreak

For all these cases the fundamental chamber $\M$ is a closed
polygon on the hyperbolic plane with angles respectively:
\newline
$A_{1,0}:\ \pi/2,\,0,\,\pi/3;$\ \
$A_{1,I}:\ 0,\,\pi/3,\,\pi/3;$\ \
$A_{1,II}:\ 0,\,0,\,0;$\ \
$A_{1,III}:\ 0,\,\pi/2,\,0,\,\pi/2,0;$
\newline
$A_{2,0}:\ 0,\,\pi/2,\,0;$\ \
$A_{2,I}:\ 0,\,\pi/2,\,0,\,\pi/2;$\ \
$A_{2,II}:\ 0,\,0,\,0,\,0;$\ \
\newline
$A_{2,III}:\ 0,\,\pi/2,\,0,\,\pi/2,\,0,\,\pi/2,\,0,\,\pi/2;$
\newline
$A_{3,0}:\ 0,\,\pi/3,\,0;$\ \
$A_{3,I}:\ 0,\,\pi/3,\,0,\,\pi/3;$\ \
$A_{3,II}:\ 0,\,0,\,0,\,0,\,0,\,0;$\ \
\newline
$A_{3,III}: \ 0,\,\pi/3,\,0,\,\pi/3,\,0,\,\pi/3,\,0,\,\pi/3,\,0,\,
\pi/3,\,0,\,\pi/3$.
\newline
All these polygons are touching a circle with the center $\br_{++}\rho$ 
where $\rho$ is the Weyl vector. It shows the geometrical sense of the 
Weyl vector.

For 9 generalized Cartan matrices $A_{i,j}$ where $i=1,\,2,\,3$ and
$j=0,\,I,\,II$ we found automorphic forms $\Phi$ for data (4), (5),
found their product formulae \thetag{1.18},  
and thus constructed the corresponding (automorphic) 
Lorentzian Kac--Moody algebras. 
See \cite{GN1}, \cite{GN2}, \cite{GN4}---\cite{GN6}. 
It is interesting that some of these automorphic forms
were well-known. For example, the automorphic form $\Phi$ for
$A_{1,II}$ is classical. It has the weight $5$ and is the product
of all even theta-constants of genus 2 (there are 10 of them). 
It is automorphic with respect to 
$Sp_4(\bz)$ with some quadratic character and
gives the discriminant of curves of genus $2$.
The automorphic form $\Phi$
for $A_{1,0}$ has the weight $35$ and is automorphic with respect to
$Sp_4(\bz)$. This automorphic form was found by Igusa $30$ years
ago and is $Sp_4(\bz)$-automorphic form of the smallest odd weight.
For both these automorphic forms, 
we found their product expansions \thetag{1.18} which
were not known. Here we use isomorphism of IV type
symmetric domain of dimension 3 with
Siegel upper-half plane of genus 2.

All other automorphic forms $\Phi$ for generalized Cartan matrices
$A_{1,0}$---$A_{3,II}$ were not known. Let us give one of them. 

We give the automorphic form $\Phi$ for $A_{3,II}$. 
For this case $T=2U(12)\oplus \langle 2 \rangle=U(12)\oplus S$ 
where $S=U(12)\oplus \langle 2 \rangle$ gives the datum (1).  
We use bases $f_2,\,\hat{f}_3,\,f_{-2}$ of the lattice $S$ with 
the Gram matrix
$$
\left(\matrix
0&0&-12\\
0&2&0\\
-12&0&0
\endmatrix
\right)
$$
and corresponding coordinates
$z_1, z_2, z_3$ of $S\otimes \bc$. The Weyl group $W$ in (2) is generated 
by reflections in all elements with square $2$ of $S$. The set $P$ in (3) 
equals  
$$
\split
P=\{&\alpha_1=(0,1,0),\,\alpha_2=(0,-1,1),\,\alpha_3=(1,-5,2),\,
\alpha_4=(2,-7,2),\\
&\alpha_5=(2,-5,1),\,\alpha_6=(1,-1,0)\}.
\endsplit
$$
It has the Gram matrix $A_{3,II}$. The Weyl vector 
$\rho=({1\over 6},\,-{1\over 2},\,{1\over 6})$. 
The automorphic form $\Phi$ is an  
automorphic cusp form $\Delta_1$ 
with respect to the group $G=O^+(T)$ with 
some character of order $6$. It has the smallest possible weight $1$   
and has the Fourier expansion and the product expansion  
$$
\align
\Delta_1(z_1, z_2, z_3)&=
\sum_{M\ge 1}
\sum
\Sb
m >0,\,l\in \bz\\
\vspace{0.5\jot} n,\,m\equiv 1\,mod\,6\\
\vspace{0.5\jot}
4nm-3l^2=M^2
\endSb
\hskip-4pt
\biggl(\dsize\frac{-4}{l}\biggr)
\biggl(\dsize\frac{12}{M}\biggr)
\sum\Sb a|(n,l,m)\endSb \biggl(\dsize\frac{6}{a}\biggr)
q^{n/6}r^{l/2}s^{m/6}\\
{}&=
q^{1/6}r^{1/2}s^{1/6}
\prod
\Sb n,\,l,\,m\in \Bbb Z\\
\vspace{0.5\jot}
(n,l,m)>0
\endSb
\bigl(1-q^{n} r^{l} s^{m}\bigr)^{f_{3}(nm,l)}
\tag{1.20}
\endalign
$$
where
$q=\exp{(24\pi i z_1)}$,
$r=\exp{(4\pi i z_2)}$,
$s=\exp{(24\pi i z_3)}$
and
$$
\biggl(\frac{-4}{l}\biggr)=\cases \pm 1 &\text{if }
l\equiv \pm 1\ \hbox{mod}\ 4\\
\hphantom{\pm}0 &\text{if }
l\equiv \hphantom{\pm} 0 \ \hbox{mod}\ 2
\endcases,\ \
\biggl(\frac{12}M\biggl)=
\cases
\hphantom{-}1\  \text{ if }\  M\equiv \pm 1 \operatorname{mod} 12\\
-1\ \text{ if }\  M\equiv \pm 5 \operatorname{mod} 12\\
\hphantom{-}0\  \text{ if }\  (M,12)\ne 1
\endcases,
$$
$$
\biggl(\frac{\,6\,}{\,a\,}\biggr)=\cases \pm 1 &\text{if }
a\equiv \pm 1\ \hbox{mod}\ 6\\
\hphantom{\pm}0 &\text{if }
(a,6)\not=1\ .
\endcases
$$
The multiplicities $f_3(nm,l)$ of the infinite product are
defined by a weak Jacobi form
$\phi_{0,3}(\tau,z)=\dsize\sum_{n\ge 0,\,l\in \bz}f_3(n,l)q^nr^l$
of weight $0$ and index $3$ with integral Fourier coefficients:
$$
\phi_{0,3}(\tau ,z)  =   r^{-1}
\biggl(\prod_{n\ge 1}(1+q^{n-1}r)(1+q^{n}r^{-1})(1-q^{2n-1}r^2)
(1-q^{2n-1}r^{-2})\biggr)^2
\tag{1.21}
$$
where $q=\exp{(2\pi i \tau )}$, $\text{im}~\tau>0$, and 
$r=\exp{(2\pi i z)}$. The divisor of $\Delta_1$ is the 
sum with multiplicities one
of all quadratic divisors orthogonal to elements of
$T$ with square $2$. The $S$, $W$, $P$ and $\Delta_1$  
define the generalized Lorentzian Kac--Moody superalgebra $\geg$ 
with the denominator identity \thetag{1.20}. 

To construct the automorphic form $\Delta_1$, we use the {\it arithmetic 
lifting of Jacobi forms on IV type symmetric domains} constructed in 
\cite{G1}---\cite{G4} 
(it gives the sum part of \thetag{1.20}). To find the product part of  
\thetag{1.20},  we use the {\it Borcherds lifting}  
\cite{B5} which is the exponential analog of the arithmetic 
lifting. See details in \cite{GN6}.

\subhead
1.5. Physical applications
\endsubhead
Lorentzian Kac--Moody algebras we have considered above and
corresponding automorphic forms $\Phi$ found some very interesting
applications in Physics: String Theory, Mirror Symmetry and others. 
We refere a reader to a nice review \cite{M} and references there.  
E. g. see \cite{Ca}, \cite{CCL}, 
\cite{DVV}, \cite{HM1}---\cite{HM2}, \cite{Ka1}---\cite{Ka2}.  
Roughly speaking the Lorentzian Kac--Moody algebras are related with 
symmetries of Fundamental Physical Theories.

\subhead
1.6. An interesting problem
\endsubhead

In data (1)---(5) above, existence of the Weyl vector
$\rho\in S\otimes \bq$ (or the lattice Weyl vector)
is very important. It is equivalent considering automorphic forms
$\Phi$ on IV type Hermitian symmetric domains. It would be
interesting to extend the Theory of Lorentzian Kac--Moody algebras
above to cases when the lattice Weyl vector $\rho$ does not exist.
It seems, for this more general case, one should consider
automorphic forms in some more general sense. On the other
hand, this more general theory will lose some finiteness properties.  
It would be a pity.

\enddocument

\end